\documentclass[11pt, a4paper]{article}
\usepackage{amsmath,amsthm,amssymb,amsfonts,esint}
\usepackage{a4wide}
\newtheorem{theorem}{Theorem}
\newtheorem{lemma}[theorem]{Lemma}
\newtheorem{proposition}[theorem]{Proposition}

\theoremstyle{remark}

\newtheorem{remark}[theorem]{Remark}

\def\Xint#1{\mathchoice
{\XXint\displaystyle\textstyle{#1}}%
{\XXint\textstyle\scriptstyle{#1}}%
{\XXint\scriptstyle\scriptscriptstyle{#1}}%
{\XXint\scriptscriptstyle\scriptscriptstyle{#1}}%
\!\int}
\def\XXint#1#2#3{{\setbox0=\hbox{$#1{#2#3}{\int}$ }
\vcenter{\hbox{$#2#3$ }}\kern-.6\wd0}}

\def\dashint{\Xint-}

\usepackage{graphicx}
\usepackage[usenames,dvipsnames]{color}
\usepackage[colorlinks=true, pdfstartview=FitV, linkcolor=black, citecolor=blue, urlcolor=blue]{hyperref}


\newcommand{\R}{{\mathbb R}}

\newcommand{\mB}{\mathcal{B}}
\newcommand{\mG}{\mathcal{G}}

\newcommand{\mM}{\mathcal{M}}
\newcommand{\mO}{\mathcal{O}}
\newcommand{\mU}{\mathcal{U}}
\newcommand{\mV}{\mathcal{V}}
\newcommand{\mW}{\mathcal{W}}

\newcommand{\pa}{{\partial}}
\newcommand{\na}{{\nabla}}

\newcommand{\eps}{{\varepsilon}}

\def\curl{\hbox{curl \!}}
\def\div{\hbox{div \!}}

\def\mspace{\medskip \noindent}

\newcommand{\dd}{\, \mathrm{d}}

\newcommand{\E}{{\mathbb E}}

\definecolor{darkgreen}{rgb}{0,0.5,0}

\title{Mild assumptions for the derivation  of  \\ Einstein's effective viscosity formula} 
\author{David G\'erard-Varet, Richard M. H\"ofer}

\begin{document}
\maketitle

\begin{abstract}
We provide a rigorous derivation of Einstein's formula for the effective viscosity of dilute suspensions of $n$ rigid balls, $n \gg 1$, set in a volume of size $1$. So far, most justifications were carried under a strong assumption on the minimal distance between the balls: $d_{min} \ge c n^{-\frac{1}{3}}$, $c > 0$. We relax this assumption into a set of two much weaker conditions: one expresses essentially that the balls do not overlap, while the other one gives a control of the number of balls  that are close to one another. In particular, our analysis covers  the case of suspensions modelled by standard Poisson processes with almost minimal hardcore condition.  
\end{abstract}

\section{Introduction}
Mixtures of particles and fluids, called {\em suspensions}, are involved in many natural phenomena and industrial processes. The understanding of their rheology, notably the so-called {\em effective viscosity} $\mu_{eff}$ induced by the particles, is therefore crucial. Many experiments or simulations have been carried out to determine $\mu_{eff}$ \cite{Guaz}. For $\lambda$ large enough, they seem to exhibit some generic behaviour, in terms of the ratio between the solid volume fraction $\lambda$ and the  maximal flowable solid volume fraction $\lambda_c$, {\it cf.} \cite{Guaz}. Still, a theoretical derivation of the relation $\mu_{eff} = \mu_{eff}(\lambda/\lambda_c)$ observed experimentally is missing, due to the complex interactions involved: hydrodynamic interactions, direct contacts, \dots Mathematical works related to the analysis of suspensions are mostly limited to the {\em dilute regime}, that is when $\lambda$ is small. 

\mspace
In these mathematical works, the typical model under consideration is as follows. One considers $n$ rigid  balls $B_i = \overline{B(x_i, r_n)}$, $1 \le i \le n$,  in a fixed compact subset of  $\R^3$, surrounded by a viscous fluid.  

The inertia of the fluid is neglected, leading to the Stokes equations
\begin{equation}
\label{Sto}
\left\{
\begin{aligned}
-\mu \Delta u_n + \na  p_n & = f_n, \quad x \in \Omega_n = \R^3 \setminus \cup B_i , \\ 
\div u_n & = 0, \quad x \in \Omega_n , \\
u_n\vert_{B_i} & = u_{n,i}  + \omega_{n,i} \times (x-x_i). 
\end{aligned}
\right.
\end{equation}
The last condition expresses a no-slip condition at the  rigid spheres, where  the velocity is given by some translation velocities $u_{n,i}$ and some rotation vectors $\omega_{n,i}$, $1 \le i  \le n$. We neglect the inertia of the balls:  the $2n$ vectors $u_{n,i}, \omega_{n,i}$ can then be seen  as Lagrange multipliers for the $2n$ conditions 
\begin{equation}
\label{Sto2}
\begin{aligned}
\int_{\pa B_i} \sigma_\mu(u,p) \nu & = - \int_{B_i} f_n , \quad \int_{\pa B_i} \sigma_\mu(u,p) \nu \times (x-x_i)  = - \int_{B_i} (x-x_i) \times f_n 
\end{aligned}
\end{equation}
where $\sigma_\mu = 2\mu D(u) \nu - p \nu$ is the usual Newtonian tensor, and $\nu$ the normal vector pointing outward $B_i$.

\mspace
 The general belief is that one should be able to replace \eqref{Sto}-\eqref{Sto2} by an effective Stokes model, with a modified viscosity taking into account the average effect of the particles:  
 \begin{equation}
\label{Stoeff}
\left\{
\begin{aligned}
-\div (2 \mu_{eff}  D u_{eff} ) + \na  p_{eff} & = f, \quad x \in \R^3, \\
\div u_{eff} & = 0, \quad x \in \R^3,
\end{aligned}
\right.
\end{equation}
with $D = \frac{1}{2}(\na + \na^t)$ the symmetric gradient. Of course, such average model can only be obtained asymptotically, namely when the number of particles $n$ gets very large. Moreover, for averaging to hold, it is very natural to impose some averaging on the distribution of the balls itself. Our basic hypothesis will therefore be the existence of a limit density, through 
\begin{equation} \label{A0} \tag{A0}
\frac{1}{n} \sum_i \delta_{x_i} \xrightarrow[n \rightarrow +\infty]{} \rho(x) dx \quad \text{weakly in the sense of measures}
\end{equation}
where $\rho \in L^\infty(\R^3)$ is assumed to be zero outside a smooth open bounded set $\mO$. After playing on the length scale, we can always assume that $|\mO| = 1$. Of course,  we expect $\mu_{eff}$ to be different from $\mu$ only in this region $\mO$ where the particles are present. 

\mspace
 The volume fraction of the balls is then given by $\lambda = \frac{4\pi}{3} n r_n^3$. We shall consider the case where $\lambda$ is small (dilute suspension), but independent of $n$ so as to derive a non-trivial effect as $n \rightarrow +\infty$. The mathematical questions that follow are: 
 \begin{itemize}
 \item Q1 :  Can we approximate system \eqref{Sto}-\eqref{Sto2} by a system of the form \eqref{Stoeff} for large $n$? 
 \item Q2 : If so, can we provide a formula for $\mu_{eff}$ inside $\mO$? In particular, for small $\lambda$, can we derive an expansion 
 $$ \mu_{eff} = \mu + \lambda \mu_1 + \dots  \quad ? $$ 
 \end{itemize}
 Regarding Q1, the only work we are aware of is the recent paper \cite{DuerinckxGloria19}. It shows that $u_n$ converges to the solution $u_{eff}$ of an effective model of the type \eqref{Sto2}, under two natural conditions: 
 \begin{enumerate}
 \item[i)]  the balls satisfy the separation condition  $\inf_{i \neq j} |x_i - x_j| \ge M \ r_n$, $M  > 2$. Note that this is a slight reinforcement of the natural constraint that the balls do not overlap. 
 \item[ii)] the centers of the balls are obtained from a stationary ergodic point process. 
 \end{enumerate}
 We refer to \cite{DuerinckxGloria19} for all details. Note that in the scalar case, with the Laplacian instead of the Stokes operator, similar results can be found in \cite[paragraph 8.6]{MR1329546}. 
 
 \mspace
 Q2, and more broadly quantitative aspects of dilute suspensions, have been studied for long. The pioneering work is due to Einstein \cite{Ein}. By {\em neglecting the interaction between the particles}, he computed a first order approximation of the effective viscosity of homogeneous suspensions: 
$$ \mu_{eff} =  (1 +  \frac{5}{2} \lambda) \mu \quad \text{ in } \mO.$$
This celebrated  formula was confirmed experimentally afterwards.  It was later extended to the inhomogenous case,  with formula 
\begin{equation}  \label{Almog-Brenner}
\mu_{eff}   = (1 +  \frac{5}{2} \lambda \rho) \mu, 
\end{equation}
see \cite[page 16]{AlBr}.  Further works investigated the $O(\lambda^2)$ approximation of the effective viscosity, {\it cf.} \cite{BaGr1}
and the recent analysis \cite{DGV_MH, GerMec20}. 

\mspace
Our concern in the present paper is the justification of Einstein's formula. To our knowledge, the  first rigorous studies on this topic are \cite{MR813656} and  \cite{MR813657}: they rely on homogenization techniques, and are restricted to suspensions that are periodically distributed in a bounded domain. A more complete  justification, still in the periodic setting but based on variational principles,  can be found in  \cite{MR2982744}. Recently,  the periodicity assumption was relaxed in \cite{HiWu}, \cite{NiSc}, and replaced by an assumption on the minimal distance: 
\begin{equation} \label{A1} \tag{A1}
\text{There exists an absolute constant $c$, such that } \quad \forall n, \forall 1 \le i \neq j \le n, \quad |x_i - x_j| \ge c n^{-\frac{1}{3}}. 
\end{equation}
For instance, introducing the solution $u_{E}$ of the Einstein's approximate model 
\begin{equation} \label{Sto_E}
-\div (2 \mu_E Du_E) + \na p_E = f, \quad \div u = 0  \quad \text{ in } \:  \R^3 
\end{equation}
with $\mu_E =  (1 +  \frac{5}{2} \lambda \rho) \mu$, it is shown in \cite{HiWu} that for all $ 1 \le p < \frac{3}{2}$,  
$$ \limsup_{n \to \infty} ||u_n - u_E||_{L^p_{loc}(\R^3)}  = O(\lambda^{1+\theta}), \quad \theta = \frac{1}{p} - \frac{2}{3}. $$
We refer to \cite{HiWu} for refined statements, including quantitative convergence in $n$ and treatment of polydisperse suspensions. 

\mspace
Although it is a substantial gain over the periodicity assumption, hypothesis \eqref{A1} on the minimal distance  is still strong. In particular, it is much more stringent that the condition that the rigid balls can not overlap. Indeed, this latter condition reads: $\forall i \neq j$, $|x_i - x_j| \ge 2 r_n$, or equivalently $|x_i - x_j| \ge c \, \lambda^{1/3} n^{-\frac{1}{3}}$, with $c = 2 (\frac{3\pi}{4})^{1/3}$. It follows from \eqref{A1} at small $\lambda$. On the other hand, one could argue that  a simple non-overlapping condition is not enough to ensure the validity of Einstein's formula. Indeed, it is based on neglecting  interaction between particles, which is incompatible with too much clustering in the suspension. Still, one can hope that if the balls are not too close from one another {\em on average}, the formula still holds.

\mspace
This is the kind of result that we prove here.  Namely, we shall replace \eqref{A1} by a set of two relaxed conditions: 
\begin{align}
\label{B1} \tag{B1} 
& \text{There exists $M >2$, such that}  \quad \forall n, \: \forall 1 \le i \neq j \le n, \quad |x_i - x_j| \ge M r_n. \\
 \label{B2} \tag{B2}
& \text{There exist $C,\alpha > 0$, such that} \quad  \forall \eta > 0,  \quad  \#\{i, \:  \exists j, \:  |x_i - x_j| \le \eta n^{-\frac13}\} \le C \eta^\alpha n 
\end{align}
Note that  \eqref{B1} is slightly stronger than the non-overlapping condition, and was already present in the work \cite{DuerinckxGloria19}  to ensure the existence of an effective model. It is possible to relax this condition into a moment bound on the particle separation, see Remark \ref{rem:recentbibli} and Section  \ref{sec:B1}. 
As regards \eqref{B2}, one can show that it is satisfied almost surely as $n \to \infty$ in the case when the particle positions 
are generated by a stationary ergodic point process if the process does not favor too much close pairs of points. In particular, it is satisfied by
a (hard-core) Poisson point process for $\alpha = 3$. 
Moreover, \eqref{B2} is satisfied for $\alpha = 3$ with probability tending to $1$ as $n \to \infty$ for independent and identically distributed particles. We postpone further discussion to Section \ref{sec:prob}.

\mspace
Under these general assumptions, we obtain: 
\begin{theorem} \label{main}
Let  $\lambda > 0$, $f \in L^1(\R^3) \cap L^\infty(\R^3)$. For all $n$, let $r_n$ such that $\displaystyle \lambda = \frac{4\pi}{3} n r_n^3$,  let $f_n \in L^{\frac65}(\R^3)$,  and $u_n$ in $\displaystyle \dot{H}^1(\R^3)$ the solution of \eqref{Sto}-\eqref{Sto2}.  Assume \eqref{A0}-\eqref{B1}-\eqref{B2}, and that $f_n \rightarrow f$ in $L^{\frac65}(\R^3)$. Then, there exists $p_{min} > 1$ such that for any $p < p_{min}$,  any $q <  \frac{3 p_{min}}{3 - p_{min}}$, one can find $\delta > 0$ with the estimate 
 $$ ||\na (u - u_E)||_{L^p(\R^3)} + \limsup_{n \rightarrow +\infty} ||u_n - u_E||_{L^q(K)} = O(\lambda^{1+\delta}), \quad \forall K \Subset \R^3, \quad \text{as } \: \lambda \rightarrow 0, $$
where   $u$  is any weak accumulation point of $u_n$ in  $\displaystyle \dot{H}^1(\R^3)$ and 
 $u_E$ satisfies Einstein's approximate model \eqref{Sto_E}. 
 \end{theorem}
\noindent Here, we use the notation $\dot H^1(\R^3)$ for the homogeneous Sobolev space 
	$\dot H^1(\R^3) = \{ w \in L^6(\R^3) : \nabla w \in L^2(\R^3)\}$ equipped with the $L^2$ norm of the gradient.
	
\begin{remark} \label{rem:exponents}
	The following explicit formula for $p_{min}$ and $\delta$ will be obtained in the proof of the theorem: 
	\begin{align*}
		p_{min} = 1 + \frac{\alpha}{6 + \alpha}, \qquad \delta = \frac 1 r - \frac{6}{6 + (2-r)\alpha} \qquad r = \max\left\{p,\frac{3q}{3+q} \right\}.
	\end{align*}
\end{remark}

\begin{remark} \label{rem:recentbibli}
	Since the preprint of our paper, several further results have appeared which we briefly discuss in this remark.

	In \cite[version 1]{DuerinckxGloria20}, an extensive study of the effective viscosity at low volume fraction was performed in the context of  stationary ergodic particle configurations, under suitable versions of \eqref{B1}-\eqref{B2}. It includes results on  the $O(\lambda^2)$ and higher order corrections, see also the recent paper \cite{Gerard-Varet20}. As regards the $O(\lambda)$ Einstein's formula, a result analogous to Theorem \ref{main} was shown with methods of a more probabilistic flavour. 
	
	It was subsequently shown in \cite{Duerinckx20} and \cite[version 2]{DuerinckxGloria20} that both the existence of an effective viscosity and the Einstein's formula hold when relaxing condition \eqref{B1} into a moment bound on the particle separation. We will argue in Section \ref{sec:B1} that our main result still holds under similar milder assumption.
	
	Finally, in \cite{HoeferSchubert20}, results have been obtained concerning the coupling of Einstein's formula to the time evolution of sedimenting particles.
\end{remark}

\mspace
The rest of the paper is dedicated to the proof of Theorem \ref{main}.  

\section{Main steps of proof}
To prove Theorem \ref{main}, we shall rely on an enhancement of  the general strategy explained  in  \cite{DGV}, to justify  various effective models for conducting and fluid media. Let us point out that one of the examples considered in \cite{DGV} is the scalar version of \eqref{Sto}-\eqref{Sto2}. It leads to a proof of a scalar analogue of Einstein's formula, under assumptions \eqref{A0}, \eqref{B1}, plus  an abstract assumption intermediate between \eqref{A1} and \eqref{B2}. We refer to the discussion at the end of \cite{DGV} for more details. Nevertheless, to justify the effective fluid model \eqref{Sto_E} under  the mild assumption \eqref{B2} will require several new steps. The main difficulty will be to handle particles that are close to one another, and will involve sharp $L^p$ estimates similar to those of \cite{GerMec20}.

\mspace
Concretely, let $\varphi$ be a smooth and compactly supported divergence-free vector field. For each $n$, we introduce the solution $\phi_n \in \dot{H}^1(\R^3)$ of 
\begin{equation} \label{Sto_phi}
\begin{aligned}
- \div(2\mu D \phi_n) + \na q_n & =   \div (5 \lambda \mu \rho D \varphi) \: \text{ in } \: \Omega_n, \\
 \div \phi_n & = 0  \: \text{ in } \: \Omega_n, \\
 \phi_n  & = \varphi + \phi_{n,i} + w_{n,i} \times (x-x_i)  \: \text{ in } \:  B_i, \: 1 \le i \le n
\end{aligned} 
\end{equation}
where the constant vectors $\phi_{n,i}$, $w_{n,i}$ are associated to the constraints 
\begin{equation} \label{Sto2_phi}
\begin{aligned}
\int_{\pa B_i} \sigma_\mu(\phi_n,q_n) \nu  & = - \int_{\pa B_i} 5 \lambda \mu \rho D \varphi \nu, \\
 \int_{\pa B_i}(x-x_i) \times  \sigma_\mu(\phi_n,q_n) \nu & = - \int_{\pa B_i}(x-x_i) \times  5 \lambda \mu \rho D \varphi \nu.
 \end{aligned}
\end{equation}
Testing \eqref{Sto} with  $\varphi - \phi_n$, we find  after a few integration by parts that 
$$
\int_{\R^3} 2\mu_E Du_n : D \varphi = \int_{\R^3}  f_n \cdot \varphi - \int_{\R^3}  f_n \cdot \phi_n.
$$
Testing \eqref{Sto_E} with $\varphi$, we find 
$$
\int_{\R^3} 2\mu_E Du_E  :  D \varphi = \int_{\R^3}  f \cdot \varphi. 
$$
Combining both, we end up with 
\begin{equation} \label{weak_estimate}
\int_{\R^3} 2\mu_E D(u_n - u_E) :  D \varphi =  \int_{\R^3}  (f_n - f) \cdot \varphi  - \int_{\R^3} f_n \cdot \phi_n. 
\end{equation}
We remind that vector fields $u_n, u_E, \phi_n$  depend implicitly on $\lambda$.

\mspace
 The main point will be to show 
\begin{proposition} \label{main_prop}
There exists $p_{min} > 1$ such that for all $p < p_{min}$, there exists $\delta > 0$ and $C > 0$,   independent of $\varphi$,  such that 
\begin{equation} \label{estimateR}
 \limsup_{n \to \infty}  \big| \int_{\R^3} f_n \cdot \phi_n \big|  \le C \lambda^{1+\delta} ||\na \varphi||_{L^{p'}}, \quad p' = \frac{p}{p-1}. 
\end{equation}
 \end{proposition}
\noindent
Let us show how the theorem follows from the proposition. First, by standard energy estimates, we find that $u_n$ is bounded in $\dot{H}^1(\R^3)$ uniformly in $n$.  Let $u = \lim u_{n_k}$ be a weak accumulation point of $u_n$ in this space. Taking the limit  in \eqref{weak_estimate}, we get 
$$ \int_{\R^3} 2\mu_E D(u - u_E) :  D \varphi =  \langle R, \varphi  \rangle $$
where  $\langle R , \varphi \rangle =   \lim_{k \rightarrow +\infty} \int_{\R^3} f_{n_k} \cdot \phi_{n_k}$. 
Recall that $\varphi$ is an arbitrary smooth and compactly supported divergence-free  vector field and that such functions are dense in the homogeneous Sobolev space of divergence-free functions $\dot{W}^{1,p}_\sigma$.
Thus, Proposition \ref{main_prop} implies that $R$ is an element of  $\dot{W}_\sigma^{-1,p}$ with $||R||_{\dot{W}_\sigma^{-1,p}} = O(\lambda^{1+\delta})$. Moreover, the previous identity is the weak formulation of 
$$ - \div(2 \mu_E D(u - u_E)) + \na q = R, \quad \div (u - u_E) = 0 \quad \text{ in } \: \R^3.    $$
Writing these Stokes equations with non-constant viscosity as 
$$ -\mu \Delta (u - u_E) + \na q = R  + \div (5 \lambda \mu \rho D(u-u_E)),  \quad \div (u - u_E) = 0 \quad \text{ in } \: \R^3.    $$
and using standard estimates for this system, we get 
$$ ||\na (u-u_E)||_{L^p} \le C \left(||R||_{\dot{W}^{-1,p}_\sigma} +  \lambda ||\na (u-u_E)||_{L^p} \right). $$
For $\lambda$ small enough, the last term is absorbed by the left-hand side, and finally 
$$ ||\na (u-u_E)||_{L^p(\R^3)} \le  C\lambda^{1+\delta}$$
which implies the first estimate of the theorem. Then, by Sobolev imbedding, for any $q \le \frac{3p}{3-p}$, and any compact $K$, 
\begin{equation} \label{Sob_imbed}
 ||u-u_E||_{L^q(K)} \le C_{K,q}  \, \lambda^{1+\delta}.
 \end{equation} 
We claim  that 
$ \limsup_{n \to \infty} ||u_n-u_E||_{L^q(K)} \le C_{K,q}  \, \lambda^{1+\delta}$. Otherwise, there exists a subsequence $u_{n_k}$ and  $\eps > 0$  such that $\displaystyle ||u_{n_k} - u_E||_{L^q(K)} \ge  C_{K,q} \,  \lambda^{1+\delta} + \eps$ for all $k$. Denoting by $u$ a (weak) accumulation point of $u_{n_k}$ in $\dot{H}^1$, Rellich's theorem implies that, for a subsequence still denoted $u_{n_k}$, $||u_{n_k} - u||_{L^q(K)} \rightarrow 0$, because $q < 6$ (for $p_{min}$ taken small enough).  Combining this  with \eqref{Sob_imbed}, we reach a contradiction. As $p$ is arbitrary in $(1, p_{min})$, $q \leq \frac{3p}{3-p}$ is arbitrary in $(1,  \frac{3 p_{min}}{3 - p_{min}})$. The last estimate of the theorem is proved.

\mspace
It remains to prove Proposition \ref{main_prop}. 
Therefore, we need a better understanding of the solution $\phi_n$ of \eqref{Sto_phi}-\eqref{Sto2_phi}. Neglecting any interaction between the balls, a natural attempt is to approximate   $\phi_n$ by
\begin{equation} \label{approx_phi}
\phi_n \approx  \phi_{\R^3} + \sum_i \phi_{i,n}
\end{equation}
where $\phi_{\R^3}$ is the solution of 
\begin{equation} \label{eq_phi_R3}
 - \mu \Delta \phi_{\R^3} + \na p_{\R^3} =  \div(5 \lambda \mu \rho D\varphi), \quad \div \phi_{\R^3}  = 0 \quad \text{in } \: \R^3 
 \end{equation}
and $\phi_{i,n}$ solves 
\begin{equation} \label{eq_phi_i_n}
 - \mu \Delta \phi_{i,n} + \na p_{i,n} = 0, \quad \div\phi_{i,n} = 0  \quad \text{ outside } \: B_i,  \quad \phi_{i,n}\vert_{B_i}(x) =  D \varphi(x_i) \, (x-x_i) 
 \end{equation}
Roughly, the idea of  approximation \eqref{approx_phi} is that $\phi_{\R^3}$ adjusts to the source term in \eqref{Sto_phi}, while for all $i$,  $\phi_{i,n}$ adjusts to the boundary condition  at the ball $B_i$. Indeed, using a Taylor expansion of $\varphi$ at $x_i$, and splitting $\na \varphi(x_i)$ between its symmetric and skew-symmetric part,  we find
$$ \phi_n\vert_{B_i}(x)  \approx D\varphi(x_i) \, (x- x_i)  \: + \:  \text{\em rigid vector field} =  \phi_{i,n}\vert_{B_i}(x) \: + \:  \text{\em rigid vector field}. $$
Moreover, $\phi_{i,n}$ can be shown to generate no force and torque, so that the extra rigid vector fields (whose role is to ensure the no-force and no-torque conditions), should be small. 

\mspace
Still, approximation \eqref{approx_phi} may be too crude :  the vector fields $\phi_{j,n}$, $j \neq i$, have a non-trivial contribution at $B_i$, and for the balls $B_j$ close to $B_i$, which are not excluded by our relaxed assumption \eqref{B1}, these contributions may be relatively big. We shall therefore modify the approximation, restricting the sum in \eqref{approx_phi} to balls far enough from the others.

\mspace
Therefore, for $\eta > 0$, we introduce a {\em good}  and a {\em bad} set of indices:  
\begin{equation} \label{good_bad_sets}
\mG_\eta = \{ 1 \le i \le n, \: \forall j \neq i, |x_i - x_j| \ge \eta n^{-\frac{1}{3}} \}, \quad \mB_\eta = \{1, \dots n\} \setminus \mG_\eta. 
\end{equation}
The good set $\mG_\eta$ corresponds to balls that are at least $\eta n^{-\frac{1}{3}}$ away from all the others. The parameter $\eta > 0$ will be specified later: we shall consider $\eta = \lambda^\theta$ for some appropriate power $0 < \theta < 1/3$.  We set 
\begin{equation} \label{def_phi_app} 
\phi_{app,n} = \phi_{\R^3}  + \sum_{i \in \mG_\eta} \phi_{i,n} 
\end{equation}
Note that $\phi_{\R^3}$ and $\phi_{i,n}$ are explicit:   
$$ \phi_{\R^3} = \mU \star \div(5 \lambda \rho D\varphi), \quad  \mU(x) = \frac{1}{8\pi} \left( \frac{I}{|x|} + \frac{x \otimes x}{|x|^3} \right)$$
and 
\begin{equation} \label{def_phi_in}
 \phi_{i,n} = r_n V[D \varphi(x_i)]\left(\frac{x-x_i}{r_n}\right) 
 \end{equation}
where for all trace-free symmetric  matrix $S$, $V[S]$ solves
$$ -\Delta V[S] + \na P[S] = 0, \: \div  V[S] = 0 \quad \text{outside } \: B(0,1), \quad  V[S](x) = Sx, \: x \in B(0,1). $$
with expressions
$$ V[S] =   \frac{5}{2} S : (x \otimes x) \frac{x}{|x|^5} + Sx \frac{1}{|x|^5}  - \frac{5}{2} (S : x \otimes x) \frac{x}{|x|^7}, \quad P[S] = 5 \frac{S : x \otimes x}{|x|^5}. $$
 Eventually, we denote 
 $$ \psi_n = \phi_n - \phi_{app,n}. $$
Tedious but straightforward calculations show that 
$$ - \div (\sigma_\mu(V[S], P[S])) = 5 \mu  S x s^1 = - \div (5 \mu S 1_{B(0,1)}) \quad \text{in } \:   \R^3  $$
where $s^1$ denotes the surface measure at the unit sphere. It follows that 
\begin{equation}
 - \mu \Delta \phi_{app,n} + \na p_{app,n} =  \div \Big( 5 \lambda \mu \rho D\varphi -  \sum_{i \in \mG_\eta} 5 \mu D \varphi(x_i) 1_{B_i} \Big), \quad \div \phi_{app,n} = 0 \quad \text{in} \:  \R^3, 
\end{equation}
Moreover, for all $1 \le i \le n$,  
\begin{align*} 
\int_{\pa B_i} \sigma_\mu(\phi_{app,n}, p_{app,n}) \nu & = - \int_{\pa B_i} 5 \lambda \mu \rho D\varphi \nu, \\
\int_{\pa B_i} (x-x_i) \times \sigma_\mu(\phi_{app,n}, p_{app,n}) \nu & = - \int_{\pa B_i} (x-x_i) \times 5 \lambda \mu \rho D\varphi \nu.
\end{align*}
Hence, the remainder $\psi_n$ satisfies 
\begin{equation} \label{Sto_psi}
\begin{aligned}
- \mu \Delta \psi_n + \na q_n & =  0 \: \text{ in } \: \Omega_n, \\
 \div \psi_n & = 0  \: \text{ in } \: \Omega_n, \\
 \psi_n  & = \varphi - \phi_{app,n} + \psi_{n,i} + w_{n,i} \times (x-x_i)  \: \text{ in } \:  B_i, \: 1 \le i \le n
\end{aligned} 
\end{equation}
where the constant vectors $\psi_{n,i}$, $w_{n,i}$ are associated to the constraints 
\begin{equation} \label{Sto2_psi}
\begin{aligned}
\int_{\pa B_i} \sigma_\mu(\psi_n,q_n) \nu  & = 0, \\
 \int_{\pa B_i}(x-x_i) \times  \sigma_\mu(\psi_n,q_n) \nu & = 0.
 \end{aligned}
\end{equation}
Estimates on $\phi_{app,n}$ and $\psi_n$ will be postponed to sections \ref{sec_app} and \ref{sec_rem} respectively. Regarding $\phi_{app,n}$, we shall prove  
\begin{proposition} \label{prop_phi_app}
For all $p \ge  1$, 
 \begin{equation} \label{estimate_phi_app}
 \limsup_{n \to \infty}  \Big|\int_{\R^3} f \cdot \phi_{app,n} \Big| \le C_{p,f} (\lambda \eta^\alpha)^{\frac{1}{p}} ||\na  \varphi||_{L^{p'}}.
 \end{equation}
 \end{proposition}
\noindent 
Regarding the remainder $\psi_n$, we shall prove
\begin{proposition}  \label{prop_psi} 
For all $1 < p < 2$, there exists   $c> 0$ independent of $\lambda$ such that for all $1 \ge \eta \ge c \lambda^{1/3}$,  
$$  \limsup_{n \to \infty}  \Big|\int_{\R^3} f \cdot \psi_n \Big|  \le C_{p,f} \lambda^{\frac 1 2} \Big(\lambda^{1+ \frac{2-p}{2p}}  \, \eta^{-\frac{3}{p}}  +  \big( \eta^\alpha \lambda \big)^{\frac{2-p}{2p}}  \Big) ||\na  \varphi||_{L^{p'}}.$$
\end{proposition}

\mspace
Let us explain how to deduce Proposition \ref{main_prop} from these two  propositions. Let $1 < p < 2$. By standard estimates, we see that $\phi_n$ is bounded uniformly in $n$ in $\dot{H}^1$. It follows that 
\begin{align*} 
 \limsup_{n \to \infty} \Big| \int_{\R^3} f_n \cdot \phi_n \Big| &=   \limsup_{n \to \infty} \Big| \int_{\R^3} f \cdot \phi_n \Big|  \le  \limsup_{n \to \infty} \Big| \int_{\R^3} f \cdot \phi_{app,n} \Big|  +   \limsup_{n \to \infty} \Big| \int_{\R^3} f \cdot \psi_n \Big| 
\\ &\le C_{p,f} \left((\lambda \eta^\alpha)^{\frac{1}{p}} + \lambda^{\frac 3 2 + \frac{2-p}{2p}} \eta^{-\frac{3}{p}}  +  \lambda^{\frac 1 2} \big( \eta^\alpha \lambda \big)^{\frac{2-p}{2p}} \right) \|\nabla \varphi \|_{L^{p'}} 
\end{align*}
To conclude, we adjust properly the parameters $p$ and $\eta$. 
We look for $\eta$ in the form $\eta = \lambda^{\theta}$, with $0 < \theta  < \frac{1}{3}$, so that the lower bound on $\eta$ needed in Proposition \ref{prop_psi} will be satisfied for small enough $\lambda$.
Then, we choose $p_{min} = 1 + \frac{\alpha}{6 + \alpha}$ and for $p < p_{min}$ we choose 
$\theta = \frac{2p}{6 + (2-p) \alpha}$. It is straightforward to check that this yields a right-hand side $\lambda^{1+\delta}$ with $\delta = \frac 1 p - \frac{6}{6 + (2-p)\alpha}$ in accordance with Remark \ref{rem:exponents}.

%

\section{Bound on the approximation} \label{sec_app}
This section is devoted to the proof  of Proposition \ref{prop_phi_app}.  We decompose 
$$\phi_{app,n} = \phi_{app,n}^1 + \phi_{app,n}^2 + \phi_{app,n}^3$$ 
where 
\begin{align*}
& - \mu \Delta \phi^1_{app,n} + \na p^1_{app,n} =  \div \Big( 5 \lambda \mu \rho D\varphi -  \sum_{1 \le i  \le n} 5 \mu D \varphi(x_i) 1_{B_i} \Big),  \quad \div \phi^1_{app,n} = 0 \quad \text{in} \:  \R^3, \\ 
& - \mu \Delta \phi^2_{app,n} + \na p^2_{app,n} =  \div \Big(   \sum_{i \in \mB_\eta}  5 \mu D \varphi(x) 1_{B_i} \Big),  \quad \div \phi^1_{app,n} = 0 \quad \text{in} \:  \R^3, \\
& - \mu \Delta \phi^3_{app,n} + \na p^3_{app,n} =  \div \Big(   \sum_{i \in \mB_\eta}  5 \mu (D \varphi(x_i) - D\varphi(x)) 1_{B_i} \Big),  \quad \div \phi^1_{app,n} = 0 \quad \text{in} \:  \R^3. 
\end{align*}
By standard energy estimates, $\phi^k_{app,n}$ is seen to be bounded in $n$ in  $\dot{H^1}$, for all $1 \le k \le 3$. We shall prove next that 
$\phi^1_{app,n}$  and $\phi^3_{app,n}$ converge in the sense of distributions to zero, while for any $f$ with $\displaystyle D ( \Delta)^{-1} \mathbb{P} f \in L^\infty$ ($\mathbb{P}$ denoting the standard Helmholtz projection),  for any $p \ge 1$, 
\begin{equation} \label{estimate_phi_app_2}
\Big|\int_{\R^3} f \cdot \phi^2_{app,n} \Big| \le C_{f,p} (\lambda \eta^\alpha)^{\frac{1}{p}}  ||\na \varphi||_{L^{p'}}, \quad p' = \frac{p}{p-1}.
\end{equation}
Proposition \ref{prop_phi_app} follows easily from those properties. 

\mspace
We start with 
\begin{lemma} 
Under assumption \eqref{A0}, $\: \sum_{1 \le i \le n} D\varphi(x_i)  \mathbf{1}_{B_i}  \: \rightharpoonup \:   \lambda \rho  D\varphi $ weakly* in $L^\infty$.   
\end{lemma}
\noindent
 {\em Proof}. As the balls are disjoint,  $|\sum_{1 \le i \le n}  D\varphi(x_i)  \mathbf{1}_{B_i}| \le ||D\varphi||_{L^\infty}$. Let $g \in C_c(\R^3)$, and denote $\delta_n = \frac{1}{n} \sum_{i} \delta_{x_i}$ the empirical measure.  We write  
 \begin{align*}
  \int_{\R^3} \sum_{1 \le i \le n} D\varphi(x_i)  \mathbf{1}_{B_i}(y) g(y) dy & =  \sum_{1 \le i \le n} D\varphi(x_i)   \int_{B(0,r_n)}  g(x_i+y) dy \\
  & = n \int_{\R^3} D\varphi(x)  \int_{B(0,r_n)}  g(x+y) dy d\delta_n(x) \\ 
  & = n r_n^3  \int_{\R^3}   \int_{B(0,1)} g(x+r_nz) dz d\delta_n(x). 
\end{align*}
The sequence of bounded continuous functions $x \rightarrow  \int_{B(0,1)} g(x+r_n z) dz$ converges uniformly to the function $x \rightarrow \frac{4\pi}{3} g(x)$ as $n \rightarrow +\infty$. We deduce:
$$ \lim_{n \to \infty}  \int_{\R^3} \sum_{1 \le i \le n}  D\varphi(x_i)   \mathbf{1}_{B_i}(y) g(y) dy =  \lim_{n \to \infty} \lambda  \int_{\R^3}   D\varphi(x)   g(x) d\delta_n(x) =   \lambda  \int_{\R^3}  D\varphi(x)   g(x) \rho(x) dx $$
where the last equality comes from \eqref{A0}.  The lemma follows by density of $C_c$ in $L^1$. 

\mspace
Let now $h  \in C^\infty_c(\R^3)$ and $v = (\Delta)^{-1} \mathbb{P} h$. We find 
\begin{align*}
\langle \phi_{app,n}^1 , h \rangle & =  \langle \phi_{app,n}^1 ,   \Delta v \rangle =  \langle   \Delta \phi_{app,n}^1 ,  v \rangle \\
& =  \int_{\R^3} \big( 5 \lambda \mu \rho D\varphi -  \sum_{1 \le i  \le n} 5 \mu D \varphi(x_i) 1_{B_i} \big) \cdot Dv \: \rightarrow 0 \quad \text{ as } \: n \rightarrow +\infty
\end{align*}
where we used the previous lemma and the fact that $Dv$ belongs to $L^1_{loc}$ and $\varphi$ has compact support.  Hence, $\phi_{app,n}^1$ converges to zero in the sense of distributions.  As regards $\phi^3_{app,n}$, we notice that 
\begin{align*}  
||\sum_{i \in \mB_\eta}  5 \mu (D \varphi(x) - D\varphi(x_i)) 1_{B_i}||_{L^1} & \le ||\na^2 \varphi||_{L^\infty}  \sum_{1 \le i \le n} \int_{B_i} |x-x_i| dx \\
& \le ||\na^2 \varphi||_{L^\infty} \lambda r_n  \rightarrow 0 \quad \text{ as } \: n \rightarrow +\infty   
\end{align*}
Using the same  duality argument as for $\phi^1_{app, n}$ (see also below), we get that $\phi^3_{app,n}$ converges to zero in the sense of distributions. 

\mspace
 It remains to show \eqref{estimate_phi_app_2}. We use a simple H\"older estimate, and write for all $p \ge 1$:  
\begin{align*}  
||\sum_{i \in \mB_\eta}  5 \mu D \varphi  1_{B_i}||_{L^1} & \le 5 \mu || \sum_{i \in \mB_\eta}  1_{B_i}||_{L^p} ||D\varphi||_{L^{p'}} = 5 \mu \big( \text{card} \ \mB_\eta  \, \frac{4\pi}{3} r_n^3 \big)^{\frac{1}{p}} ||D\varphi||_{L^{p'}} \\
& \le C (\eta^\alpha \lambda)^{\frac{1}{p}}  ||D\varphi||_{L^{p'}} 
\end{align*}
where the last inequality follows from \eqref{B2}. Denoting $v = ( \Delta)^{-1}  \mathbb{P} f$, we have this time  
\begin{align*}
\int_{\R^3} f \cdot \phi_{app,n}^2 & =   \int_{\R^3} D v \cdot \sum_{i \in \mB_\eta}  5 \mu D \varphi 1_{B_i} \le C ||Dv||_{L^\infty}  (\eta^\alpha \lambda)^{\frac{1}{p}}  ||D\varphi||_{L^{p'}} 
\end{align*}
which implies \eqref{estimate_phi_app_2}.

\section{Bound on the remainder} \label{sec_rem}
We focus here on estimates for the remainder $\psi_n = \phi_n - \phi_{app,n}$, which satisfies \eqref{Sto_psi}-\eqref{Sto2_psi}.
The proof of Proposition  \ref{prop_psi} relies on properties of the solutions of the system
\begin{equation} \label{Sto_Psi}
-\mu \Delta \psi + \na p = 0, \quad \div \psi = 0 \quad \text{ in } \: \Omega_n, \quad 
D \psi = D \tilde{\psi} \quad \text{ in } \: B_i, \quad 1 \le i \le n   
\end{equation} 
together with the constraints 
\begin{equation} \label{Sto2_Psi} 
 \int_{\pa B_i} \sigma_\mu(\psi, p)\nu =  \int_{\pa B_i} (x-x_i) \times \sigma_\mu(\psi, p)\nu = 0, \quad 1 \le i \le n.
 \end{equation}
 
More precisely, we use a duality argument to prove the following proposition, corresponding to \cite[Proposition 3.2]{GerMec20}.
\begin{proposition} \label{prop_estimate_gphi2n} 
Let $q > 3$. Then, under assumption \eqref{B1} for all  $g \in L^{q}(\R^3)$ and all $\tilde \psi \in H^1(\cup_i B_i)$, the weak solution $\psi \in \dot H^1(\R^3)$ to \eqref{Sto_Psi}-\eqref{Sto2_Psi} satisfies
\begin{align} \label{improvement}
	\left|\int_{\R^3} g \psi \right| \leq C_{g} \lambda^{\frac 1 2 } \| D \tilde{\psi} \|_{L^2(\cup B_i)}.
\end{align}
\end{proposition}
\begin{proof}
We introduce the solution $u_g$ of the Stokes equation 
\begin{equation} \label{eq_ug}
 -\Delta u_g + \na p_g = g, \quad \div g = 0, \quad \text{ in } \: \R^3.
 \end{equation}
As $g \in L^{q}$, $q>3$, $u_g \in W^{2,q}_{loc}$, so that $D(u_g)$ is continuous.    
Integrations by parts yield 
\begin{align*}
\int_{\R^3} g \psi & = \int_{\R^3}(-\Delta u_g + \na p_g) \psi = 2 \int_{\R^3} D(u_g) : D(\psi) \\
& = 2 \int_{\cup B_i} D(u_g) : D(\psi) - \sum_i \int_{\pa B_i} u_g \cdot \sigma(\psi,  p)\nu \\
& = 2 \int_{\cup B_i} D(u_g) : D(\psi) - \sum_i \int_{\pa B_i} (u_g + u^i_g + \omega^i_g \times (x-x_i)) \cdot \sigma(\psi,  p)\nu
\end{align*}
for any constant vectors $u^i_g$, $\omega^i_g$, $1 \le i \le n$, by the force-free and torque-free conditions on $\psi$.  
As $u_g + u^i_g + \omega^i_g \times (x-x_i)$ is divergence-free, one has 
$$ \int_{\pa B_i} (u_g + u^i_g + \omega^i_g \times (x-x_i)) \cdot \nu = 0.  $$
We can apply classical considerations on the Bogovskii operator: for  any $1 \le i \le n$, there exists $U_g ^i \in H^1_0(B(x_i, (M/2)r_n))$ such that 
$$ \div U_g ^i = 0 \quad \text{ in } \: B\Big(x_i, \frac{M}{2} r_n\Big), \quad U_g ^i = u_g + u^i_g + \omega^i_g \times (x-x_i)  \quad \text{ in } \: B_i $$
and with  
$$ ||\na U_g ^i||_{L^2} \le C_{i,n} ||u_g + u^i_g + \omega^i_g \times (x-x_i)||_{W^{1,2}(B_i)} $$
Furthermore, by a proper choice of $u_g ^i$ and $\omega_g^i$, we can ensure the Korn inequality:  
$$ ||u_g + u^i_g + \omega^i_g \times (x-x_i)||_{W^{1,2}(B_i)} \le c'_{i,n} ||D(u_g)||_{L^2(B_i)} $$
resulting in 
\begin{equation*} \label{control_Ug^i} 
||\na U_g ^i||_{L^2} \le C ||D(u_g)||_{L^2(B_i)} 
\end{equation*}
where the constant $C$ in the last inequality can be taken independent of $i$ and $n$ by translation and scaling arguments. Extending $U_g ^i$ by zero, and denoting $U_g = \sum U_g ^i$, we have 
\begin{equation} \label{control_Ug} 
||\na U_g||_{L^2} \le C ||D(u_g)||_{L^2(\cup B_i)} 
\end{equation}

Thus, we find 
\begin{align*}
\int_{\R^3} g \psi & = 2 \int_{\cup B_i} D(U_g) : D(\psi) - \sum_i \int_{\pa B_i} U_g \cdot \sigma(\psi,  q)\nu \\
& = 2 \int_{\R^3} D(U_g) : D(\psi) 
\end{align*}
By using \eqref{control_Ug} and Cauchy-Schwarz inequality, we end up with 
 \begin{align*}
\big| \int_{\R^3} g \psi \big| & \le C ||D(u_g)||_{L^2(\cup B_i)} \|D(\psi)\|_{L^2(\R^3)} \le C ||D(u_g)||_{L^\infty} \lambda^{\frac12}  \|D(\psi)\|_{L^2(\R^3)}
\end{align*}
Now the assertion follows from the somehow standard estimate  
\begin{equation} \label{L2_estimate_Psi} 
||\na \psi||_{L^2(\R^3)} \le C \| D \tilde{\psi} \|_{L^2(\cup B_i)} 
\end{equation} 
for a constant $C$ independent of $n$. Indeed, by  a classical variational characterization of $\psi$, we have 
$$
||\na \psi||_{L^2(\R^3)}^2 =  2 ||D \psi||_{L^2(\R^3)}^2 = \inf \big\{  2 ||D U||_{L^2(\R^3)}^2, \: D U = D \tilde{\psi} \: \text{ on } \cup_i B_i\big\}. 
$$
Thus, \eqref{L2_estimate_Psi} follows by constructing such a vector field $U$ from $\tilde \psi$ in the same manner as we constructed $U_g$ from $u_g$ above and applying \eqref{control_Ug}.
\end{proof}

\mspace
By \eqref{Sto_psi} we can apply this proposition with $g=f$, $\psi = \psi_n$ and $\tilde \psi_n = \varphi - \phi_{app,n}$.
Thus, for the proof of Proposition \ref{prop_psi}, it remains to show
\begin{equation} \label{bound.psi}
\limsup_{n \to \infty} ||D (\varphi - \phi_{app,n})||_{L^2(\cup B_i)} \le C  \Big(\lambda^{1+ \frac{2-p}{2p}}  \, \eta^{-\frac{3}{p}}  +  \big( \eta^3 \lambda \big)^{\frac{2-p}{2p}}  \Big) ||\na  \varphi||_{L^{p'}}.
\end{equation}

  We  decompose
$$\varphi - \phi_{app,n} = \tilde{\psi}^1_n + \tilde{\psi}^2_n + \tilde{\psi}^3_n $$
where  
\begin{align*}
& \forall 1 \le i \le n, \: \forall x \in B_i, \quad  \tilde{\psi}^1_n(x)  = -\phi_{\R^3}(x) - \sum_{\substack{j \neq i, \\ j \in \mG_\eta}} \phi_{j,n}(x) 
\end{align*}
and 
\begin{align*}
& \forall i \in \mG_\eta,   \: \forall x \in B_i, \quad  \tilde{\psi}^2_n(x)  = \varphi(x) -  \varphi(x_i) - \na\varphi(x_i) (x-x_i) + \Bigl(\varphi(x_i) + \frac 1 2 \curl \varphi(x_i) \times (x-x_i)\Bigr), \\
& \forall i \in \mB_\eta,  \: \forall x \in B_i, \quad  \tilde{\psi}^2_n(x)   = 0, \\
& \forall i \in \mG_\eta,   \: \forall x \in B_i, \quad  \tilde{\psi}^3_n(x) = 0, \\
&  \forall i \in \mB_\eta, \: \forall x \in B_i,  \quad \tilde{\psi}^3_n(x) = \varphi(x). 
\end{align*}
We remind that the sum in \eqref{def_phi_app} is restricted to indices $i \in \mG_\eta$ and that $\phi_{i,n}(x) = D\varphi(x_i) (x-x_i)$ for $x$ in $B_i$. This explains the distinction between $\tilde{\psi}^2_n$ and  $\tilde{\psi}^3_n$.

\mspace
The control of $\tilde \psi^2_n$ is the simplest:
\begin{align} \label{psi_2n}
 \|D \tilde \psi^2_n\|_{L^2(\cup B_i)} \le C ||D^2 \varphi||_{L^\infty} \Big( \sum_{i \in \mG_\eta} \int_{B_i} |x-x_i|^2 dx \Bigr)^{1/2}  \le C' \lambda^{1/2} r_n. 
\end{align} 
Hence, 
\begin{equation} \label{lim_psi_2n}
\lim_{n \rightarrow +\infty}  \|D \tilde \psi^2_n\|_{L^2(\cup B_i)} = 0. 
\end{equation}

Next, we estimate $\tilde \psi^3_n$.
 This term expresses the effect of the balls $\mB_\eta$ that are close to one another. By assumption  \eqref{B2}, $\text{card}\mB_\eta \le C \eta^\alpha n$.
  Thus, 
\begin{equation} \label{bound_psi_3n}
\begin{aligned}
 \|D \tilde \psi^3_n\|_{L^2(\cup B_i)}  \le C  \|1_{\cup_{i \in \mB'} B_i}\|_{L^{\frac{2p}{2-p}}(\R^3)} \|D\varphi\|_{L^{p'}(\cup_{i \in \mB'} B_i)} 
 \leq C'(\eta^\alpha \lambda)^{\frac {2-p} {2p}} \|\nabla \varphi\|_{L^{p'}}.
\end{aligned}
\end{equation}

The final step in the proof of Proposition  \ref{prop_psi} is to establish bounds on $\tilde \psi^1_n$.
We have
\begin{align} \label{bound_psi_2n}
 ||D \tilde \psi^1_n||_{L^2(\cup B_i)} & \le C \Big( \|D \phi_{\R^3}\|_{L^2(\cup B_i)}  + \Big(\sum_i  \int_{B_i} \big| \sum_{\substack{j \neq i, \\ j \in \mG_\eta}} D\phi_{j,n}\big|^2 \Big)^{1/2} \Big)  
 \end{align}
 For any $r,s < +\infty$ with $\frac{1}{r} + \frac{1}{s} = \frac{1}{2}$, we obtain 
 \begin{equation} \label{bound_phi_R3}
   \|D \phi_{\R^3}\|_{L^2(\cup B_i)} \: \le \: ||1_{\cup B_i}||_{L^r(\R^3)} ||D \phi_{\R^3}||_{L^s(\R^3)} \: \le C  ||1_{\cup B_i}||_{L^r(\R^3)} ||\lambda \rho D\varphi||_{L^s(\R^3)} 
 \end{equation}
using standard $L^s$ estimate for system \eqref{eq_phi_R3}. Hence, 
$$ \|D \phi_{\R^3}\|_{L^2(\cup B_i)} \le C' \lambda^{1+ \frac{1}{r}}  ||D\varphi||_{L^s(\R^3)}.  $$
Note that we can choose any $s >2$, this lower bound coming from the requirement $\frac{1}{r} + \frac{1}{s} = \frac{1}{2}$. Introducing $p$ such that $s= p'$, we find that for any $p < 2$, 
\begin{equation} \label{bound_Dphi_R3}
 \|D \phi_{\R^3}\|_{L^2(\cup B_i)} \le C'  \lambda^{\frac{1}{2} + \frac{1}{p}}    ||D\varphi||_{L^{p'}(\R^3)}. 
 \end{equation}
 The treatment of the second term at the r.h.s. of \eqref{bound_psi_2n} is more delicate. We write, see \eqref{def_phi_in}: 
\begin{align} \label{decompo_phi_jn}
D\phi_{j,n}(x) & = DV[D\varphi(x_j)]\Big(\frac{x-x_j}{r_n}\Big)  = \mV[D\varphi(x_j)]\Big(\frac{x-x_j}{r_n}\Big)  \:  + \:   \mW[D\varphi(x_j)]\Big(\frac{x-x_j}{r_n}\Big)
\end{align}
where $\: \mV[S] = D\Big( \frac{5}{2} S : (x \otimes x) \frac{x}{|x|^5} \Big)$, $\: \mW[S] = D \Big( \frac{Sx}{|x|^5}  - \frac{5}{2} (S : x \otimes x) \frac{x}{|x|^7} \Big)$.

\mspace
We have: 
\begin{align*}
& \sum_i  \int_{B_i} \Big| \sum_{\substack{j \neq i, \\ j \in \mG_\eta}}    \mW[D\varphi(x_j)]\Big(\frac{x-x_j}{r_n}\Big)\Big|^2  dx
\: \le \:   C \,   r_n^{10}  \, \sum_i  \int_{B_i}  \Big( \sum_{\substack{j \neq i, \\ j \in \mG_\eta}}   |D\varphi(x_j)| \,  |x-x_j|^{-5} \Big)^2  dx
\end{align*}
For all $i$, for all $j \in \mG_\eta$ with $j \neq i$, and all $(x,y) \in B_i \times B(x_j, \frac{\eta}{4} n^{-\frac{1}{3}})$,  we have for some absolute constants $c,c' > 0$: 
$$|x-x_j| \: \ge \:  c \,  |x - y|  \ge c' \,  \eta n^{-\frac{1}{3}}. $$  
Denoting $B_j^* = B(x_j,\frac{\eta}{4} n^{-\frac{1}{3}})$ We deduce 
\begin{align*}
&  \sum_i  \int_{B_i} \Big| \sum_{\substack{j \neq i, \\ j \in \mG_\eta}}    \mW[D\varphi(x_j)]\Big(\frac{x-x_j}{r_n}\Big)\Big|^2  dx \\
& \le  C \,   r_n^{10}  \sum_i \int_{B_i}  
\Big(  \sum_{\substack{j \neq i, \\ j \in \mG_\eta}}  \frac{1}{|B_j^*|} \int_{B_j^*} |x - y|^{-5} 1_{\{|x-y| > c' \eta n^{-\frac{1}{3}}\}}(x-y) |D\varphi(x_j)| dy  \Big)^2 dx\\
& \le   C' \,  n^2 \frac{r_n^{10}}{\eta^6}  \int_{\R^3} 1_{\cup B_i}(x)  \Big( \int_{\R^3} |x - y|^{-5} 1_{\{|x-y| > c' \eta n^{-\frac{1}{3}}\}}(x-y) \sum_{1 \le j \le n}  |D\varphi(x_j)| 1_{B_j^*}(y) dy  \Big)^2 dx 
\end{align*}
Using H\"older and Young's convolution inequalities, we find that for all $r,s$ with $\frac{1}{r} + \frac{1}{s} = 1$, 
\begin{align*} 
& \int_{\R^3} 1_{\cup B_i}(x)  \Big( \int_{\R^3} |x - y|^{-5}  1_{\{|x-y| > c' \eta n^{-\frac{1}{3}}\}}(x-y) \sum_{1 \le j \le n}  |D\varphi(x_j)| 1_{B_j^*}(y) dy  \Big)^2 dx  \\
& \le  ||1_{\cup B_i}||_{L^r}  \, || \big(|x|^{-5} 1_{\{|x| > c' \eta n^{-\frac{1}{3}}\}}\big) \star  \sum_{1 \le j \le n}  |D\varphi(x_j)| 1_{B_j^*} ||_{L^{2s}}^2  \\
& \le  ||1_{\cup B_i}||_{L^r}  \,  |||x|^{-5} 1_{\{|x| >  c' \eta n^{-\frac{1}{3}}\}}||_{L^1}^2 \, ||\sum_{1 \le j \le n}  |D\varphi(x_j)| 1_{B_j^*} ||_{L^{2s}}^2  \\
& \le C \lambda^{\frac{1}{r}} \,   (\eta n^{-\frac{1}{3}})^{-4} \, \Big( \sum_j |D\varphi(x_j)|^{2s} \eta^3 n^{-1} \Big)^{\frac{1}{s}}
\end{align*}
Note that, by \eqref{A0}, $\frac{1}{n}  \sum_j |D\varphi(x_j)|^t \rightarrow \int_{\R^3} |D\varphi|^t \rho$ as $n \rightarrow +\infty$. We end up with 
\begin{align*}
&  \limsup_{n \to \infty} \, \sum_i  \int_{B_i} \Big| \sum_{\substack{j \neq i, \\ j \in \mG_\eta}}    \mW[D\varphi(x_j)]\Big(\frac{x-x_j}{r_n}\Big)\Big|^2  dx 
 \le C  \, \lambda^{\frac{10}{3} + \frac{1}{r}} \, \eta^{-10+\frac{3}{s}} ||D\varphi||^2_{L^{2s}(\mO)}. 
\end{align*}
We can take any $s > 1$, which yields by setting $p$ such that  $p'=2s$: for any $p < 2$ 
\begin{equation} \label{bound_mW}
  \limsup_{n \to \infty} \, \sum_i  \int_{B_i} \Big| \sum_{\substack{j \neq i, \\ j \in \mG_\eta}}    \mW[D\varphi(x_j)]\Big(\frac{x-x_j}{r_n}\Big)\Big|^2  dx  
 \le C  \, \lambda^{\frac{10}{3} + \frac{2-p}{p}} \, \eta^{-4-\frac{6}{p}} ||D\varphi||^2_{L^{2s}(\mO)}. 
\end{equation}
To treat the first term in the decomposition \eqref{decompo_phi_jn}, we write 
$$\mV[D\varphi(x_j)]\Big(\frac{x-x_j}{r_n}\Big)  = r_n^3 \, \mM(x-x_j) \, D\varphi(x_j) $$
for $\mM$ a matrix-valued Calderon-Zygmund operator.

We use that for all $i$ and all $j \neq i$, $j \in \mG_\eta$ we have for all $(x,y) \in B_i \times B_j^\ast$ 
\begin{align*}
 |\mM(x-x_j) - \mM(x-y)| \leq C \eta n^{-1/3} |x- y|^{-4}
\end{align*}
Thus, by similar manipulations as before
\begin{align*}
&  \sum_i  \int_{B_i} \Big| \sum_{\substack{j \neq i, \\ j \in \mG_\eta}}    \mV[D\varphi(x_j)]\Big(\frac{x-x_j}{r_n}\Big)\Big|^2  dx \\
& \leq C r_n^6 \sum_i \int_{B_i}  
\Big(  \sum_{\substack{j \neq i, \\ j \in \mG_\eta}}  \frac{1}{|B_j^*|} \int_{B_j^*} \mM(x-y) 1_{\{|x-y| > c \eta n^{-\frac{1}{3}}\}}(x-y) |D\varphi(x_j)| dy  \Big)^2 dx\\
& + C \frac{\eta^2}{n^{2/3}} r_n^6 \sum_i \int_{B_i}  
\Big(  \sum_{\substack{j \neq i, \\ j \in \mG_\eta}}  \frac{1}{|B_j^*|} \int_{B_j^*} |x- y|^{-4} 1_{\{|x-y| > c \eta n^{-\frac{1}{3}}\}}(x-y) |D\varphi(x_j)| dy  \Big)^2 dx \\
& \leq  C n^2 \frac{r_n^6}{\eta^6} \,  ||1_{\cup B_i}||_{L^r}  \, || \big(\mM(x)  1_{\{|x| > \eta n^{-\frac{1}{3}}\}}\big) \star  \sum_{1 \le j \le n} \,   |D\varphi(x_j)| 1_{B_j^*} ||_{L^{2s}}^2  \\
&+C n^2 \frac{\eta^2}{n^{2/3}} \frac{r_n^6}{\eta^6} ||1_{\cup B_i}||_{L^r}  \,  |||x|^{-4} 1_{\{|x| >  c \eta n^{-\frac{1}{3}}\}}||_{L^1}^2 \, ||\sum_{1 \le j \le n}  |D\varphi(x_j)| 1_{B_j^*} ||_{L^{2s}}^2  \\
\end{align*}

As seen in \cite[Lemma 2.4]{DGV_MH}, the kernel $\mM(x)  1_{\{|x| > c \eta n^{-\frac{1}{3}}\}}$ defines a singular integral that is continuous over $L^t$ for any $1 <  t < \infty$, with operator norm bounded independently of the value $\eta n^{-\frac{1}{3}}$ (by scaling considerations). Applying this continuity property with $t=2s$, writing as before $p'=2s$, we get for all $p < 2$,
\begin{align*} 
& \limsup_{n \to \infty}   \sum_i  \int_{B_i} \Big| \sum_{\substack{j \neq i, \\ j \in \mG_\eta}}    \mV[D\varphi(x_j)]\Big(\frac{x-x_j}{r_n}\Big)\Big|^2  dx \le C \lambda^{2+ \frac{2-p}{p}} \eta^{-\frac{6}{p}}  ||D\varphi||^2_{L^{2s}(\mO)}
\end{align*}
Combining this last inequality with \eqref{decompo_phi_jn} and \eqref{bound_mW}, we finally get: for all $p < 2$, 
\begin{align} \label{bound_sum_phi_jn}
&  \limsup_{n \to \infty} \Big(\sum_i  \int_{B_i} \big| \sum_{\substack{j \neq i, \\ j \in \mG_\eta}} D\phi_{j,n}\big|^2 \Big)^{1/2}  \le C' \lambda^{1+ \frac{2-p}{2p}} \, \eta^{-\frac{3}{p}}  ||D\varphi||_{L^{p'}(\mO)}
\end{align}
Finally, if we inject \eqref{bound_Dphi_R3} and \eqref{bound_sum_phi_jn} in \eqref{bound_psi_2n}, we obtain that for any $p < 2$, 
\begin{align} \label{bound_psi_1n_final}
 \limsup_{n \to \infty} ||D \tilde \psi^1_n||_{L^2(\cup_i B_i)} & \le C \lambda^{1+ \frac{2-p}{2p}} \, \eta^{-\frac{3}{p}}   ||D\varphi||_{L^{p'}(\R^3)} 
 \end{align}
 Here we used $\eta \leq 1$.

The desired estimate \eqref{bound.psi}  follows from collecting \eqref{lim_psi_2n}, \eqref{bound_psi_1n_final} and \eqref{bound_psi_3n}. 
This concludes the proof of Proposition \ref{prop_psi}.

\section{Discussion of assumption \texorpdfstring{\eqref{B1}}{(B1)}} \label{sec:B1}

In the light of the recent paper \cite{Duerinckx20}, we will show how condition  \eqref{B1} can be replaced by the following assumption:
\begin{align} \label{B1'} \tag{B1'} 
	\forall i, \quad \rho_i := \sup_{ j \neq i} r_n^{-1} |x_i - x_j| - 2 > 0, \qquad \exists s > 1, \quad  \limsup_{n \to \infty} \frac{1}{n} \sum_i   \rho_i^{-s} < \infty.
\end{align}

%
%
We will argue that Theorem \ref{main} remains valid with $p_{\min}$ depending in addition on the power $s$ from \eqref{B1'}. More precisely, $p_{\min}$ in Remark \ref{rem:exponents} needs to be replaced by 
\begin{align} \label{p_min}
	p_{\min} = \min\left\{ 1 + \frac{\alpha}{6+\alpha}, 1 + \frac{s-1}{s+1}, \frac{3}{2}\right\}. 
\end{align} 

There are only two instances where we have used assumption \eqref{B1},  which are both contained in the proof of Proposition \ref{prop_estimate_gphi2n} :  one is to prove the estimate \eqref{L2_estimate_Psi} for the solution $\psi$ to the system
\eqref{Sto_Psi} --  \eqref{Sto2_Psi}, and the other one is to prove the analogue estimate \eqref{control_Ug}. The proof has been based on the construction of suitable functions $\Psi_i \in \dot H^1_0(B(x_i,M/2 r_n))$ with $D(\Psi_i) = D(\tilde \psi)$ in $B_i$. If we drop assumption \eqref{B1}, we can still replace the balls  $B(x_i,M/2 r_n)$ by disjoint neighborhoods $B_i^+$ satisfying the assumptions of \cite[section 3.1]{Duerinckx20} (with $I_n, I_n^+$ replaced by $B_i, B_i^{+}$). By \cite[Lemma 3.3]{Duerinckx20}, it then follows that for all $r > 2$ and all $q \ge \max(2,\frac{6r}{5r-6})$, there exists $\Psi_i \in  H^1_0(B_i^+)$, such that 
$$ \| \na \Psi_i\|_{L^2(B_i^+)} \le C_{r} \, \rho_i^{\frac{2}{r} - \frac{3}{2}} r_n^{\frac{3}{2} - \frac{3}{q}}   \| D(\tilde \psi)\|_{L^q(B_i)},  $$
Setting $\Psi = \sum_i \Psi_i$, we find that 
\begin{equation} \label{estim_naPsi}
	\|\nabla \Psi\|^2_{L^2} \leq C_{r}  \sum_i  \rho_i^{\frac{4}{r} - 3} r_n^{3-\frac{6}{q}}   \| D(\tilde \psi)\|_{L^q(B_i)}^2  \leq C_{r}  \lambda^{\frac{q-2}{q}} \left( \frac{1}{n}  \sum_i  
	\left( \rho_i^{\frac{4}{r} - 3}\right)^\frac{q}{q-2}  \right)^{\frac{q-2}{q}} \|D(\tilde \psi)\|_{L^q(\cup B_i)}^2
	\end{equation}

Note that for $s$ the exponent in \eqref{B1'},  $q > 3$ and $\frac{q}{q-2} < s$,  taking $r$ close enough to $2$, one can ensure that  $q \ge \max(2,\frac{6r}{5r-6})$ and that the first  factor at the right-hand side of \eqref{estim_naPsi} is finite. In conclusion this argument shows that Proposition \ref{prop_estimate_gphi2n} remains valid under assumption \eqref{B1'} with the estimate \eqref{improvement} replaced by 
	\begin{align} \label{improvement'}
	\left|\int_{\R^3} g \psi \right| \leq C_{g,q} \lambda^{\frac 1 2 + \frac{q - 2}{2q}}  \| D \tilde{\psi} \|_{L^{q}(\cup B_i)}.
\end{align}
It is not difficult to check that this change of the estimate still allows to conclude the argument in Section \ref{sec_rem} along the same lines as before. Indeed, whenever we used \eqref{improvement}, we also applied Hölder's estimate to replace $\| D \tilde{\psi} \|_{L^{2}(\cup B_i)}$ by a higher Lebesgue norm in order to gain powers in $\lambda$. One could say that the modified estimate \eqref{improvement'} has just partly anticipated Hölder's estimate. The additional restrictions on $q$ ($q >3$, $q < \frac{s}{s-2}$) are the reason for the additional constraints in $p_{\min}$ in \eqref{p_min}. The estimates in Section \ref{sec_rem} where we use Proposition \ref{prop_estimate_gphi2n} concern the terms $\tilde \psi_n^i$, $i=1,2,3$.
First, in the estimate for $\tilde \psi_n^2$ corresponding to \eqref{psi_2n}, we can just use \eqref{improvement'} with $q= \infty$.
Second for  $\tilde \psi_n^3$, previously estimated in \eqref{bound_psi_3n},
we use \eqref{improvement'} with $q = p'$.  
Finally, for $\tilde \psi_n^1$, if one carefully follows the estimates in Section \ref{sec_rem}, one observes that \eqref{improvement'} with $q= p'$ is again sufficient.

\section{Discussion of assumption \texorpdfstring{\eqref{B2}}{(B2)}} \label{sec:prob}

\subsection{Stationary ergodic processes}

Let $\Phi^\delta = \{y_i\}_i \subset \R^3$ be a stationary ergodic point process on $\R^3$ with intensity $\delta$ and hard-core radius $R$, i.e., $|y_i - y_j| \geq R$ for all $i \neq j$. An example of such a process is a hard-core Poisson point process, which is obtained from a Poisson point process upon deleting all points with a neighboring point closer than $R$. We refer to \cite[Section 8.1]{MR1950431} for the construction and properties of such processes.

\mspace
Assume that $\mO$ is convex and contains the origin. For $\eps > 0$, we consider the set  
\begin{align*}
 \eps \Phi^\delta \cap \mO =: \{ x^\eps_i,  i =1, \dots, n_\eps\}.
\end{align*}
Let $r < R/2$ and denote $r_\eps = \eps r$ and consider $B_i = \overline{ B(x_i,r_\eps)}$.
The volume fraction of the particles depends on $\eps$ in this case.
However, it is not difficult to generalize our result to the case when the volume fraction converges to $\lambda$ and
this holds in the setting under consideration since
\begin{align*}
	 \frac{4 \pi}{3} n_\eps r_\eps^3 \to \frac{4 \pi}{3} \delta r^3 =: \lambda(r,\delta) \quad	\text{almost surely as } \eps \to 0.
\end{align*}
Clearly, $\lambda(r,\delta) \to 0$, both if $r \to 0$ and if $\delta \to 0$.
However, the process behaves fundamentally differently in those cases. Indeed, if we take $r \to 0$ (for $\delta$ and $R$ fixed), we find that condition \eqref{A1}, which implies \eqref{B2}, is satisfied
almost surely for $\eps$ sufficiently small as
\begin{align*}
	n_\eps^{1/3} |x_i^\eps - x_j^\eps| \geq n_\eps^{1/3} \eps R \to \delta^{1/3} R.
\end{align*}

\mspace
In the case when we fix $r$ and consider $\delta \to 0$ (e.g. by randomly deleting points from a process), \eqref{A1} is in general not satisfied.
We want to characterize processes for which \eqref{B2} is still fulfilled almost surely as $\eps \to 0$.
Indeed, using again the relation between $\eps$ and $n_\eps$, it suffices to show 
\begin{align} \label{eq:B2.prob}
	\forall \eta > 0,  \quad  \#\{i, \:  \exists j, \:  |x_i - x_j| \le \eta  \eps \} \le C \eta^\alpha \delta^{1 + \frac \alpha 3} \eps^{-3}.
\end{align}
Let $\Phi^\delta_\eta$ be the process obtained from $\Phi^\delta$ by deleting those points  $y$ with $B(y,\eta) \cap \Phi^\delta = \{y\}$.
Then, the process $\Phi^\delta_\eta$ is again stationary ergodic (since deleting those points commutes with translations\footnote{In detail: let $\mathcal E_\eta$ be the operator that erases all points without a neighboring point closer than $\eta$, and let $T_x$ denote a translation by $x$.  Now, let $\mu$ be the measure for the original process $\Phi^\delta$. Then the measure for $\Phi^\delta_\eta$
is given by $\mu_\eta = \mu \circ \mathcal E_\eta^{-1}$. Since $\mathcal E_\eta T_x = T_x \mathcal E_\eta$ (for all $x$, in particular for $T_{-x} = T^{-1}_x$), 
we have for any measurable set $A$ that $T_x \mathcal E_\eta^{-1} A = \mathcal E_\eta^{-1} T_x  A$.
This immediately implies that the new process adopts stationarity and ergodicity.}), so that almost surely as $\eps \to 0$
\begin{align*}
\eps^3  \#\{i, \:  \exists j, \:  |x_i - x_j| \le \eta  \eps \} \to \E[\# \Phi^\delta_\eta \cap Q],
\end{align*}
where $Q = [0,1]^3$. Clearly, 
$$ \E[\# \Phi^\delta_\eta \cap Q]  \le  \E \sum_{y  \in \Phi^\delta \cap Q}  \sum_{y' \neq y \in \Phi^\delta} 1_{B(0,\eta)}(y' - y). $$
 We can express this expectation in terms of the 2-point correlation function $\rho^\delta_2(y,y')$ of $\Phi^\delta$ yielding
$$  \E[\# \Phi^\delta_\eta \cap Q] \le \int_{\R^6} 1_Q(y) 1_{B(0,\eta)}(y'-y) \rho^\delta_2(y,y') \dd y \dd y'. $$
Hence, \eqref{eq:B2.prob} and therefore also \eqref{B2} is in particular satisfied with $\alpha = 3$ if $\rho^\delta_2 \leq C \delta^2$
which is the case for a (hard-core) Poisson point process.

Moreover, we observe that \eqref{B2} with $\alpha < 3$ is satisfied even for processes
that favor clustering: \eqref{eq:B2.prob} holds if $\rho_2^\delta(y,y') \leq C \delta^{1 + \frac{\alpha}3} |y - y'|^{\alpha - 3}$. This means that $\rho_2^\delta$ can be quite singular at the diagonal and of much higher intesity than $\delta^2$. Examples for such clustering point processes are Neyman-Scott processes (see e.g. \cite[Section 6.3]{MR1950431}).


%

\subsection{Identically, independently distributed particles}

Focusing on assumption \eqref{B2}, we neglect the non-overlapping condition \eqref{B1} in the following, which is not satisfied for i.i.d. particles. As in the case of  hard-core Poisson point processes, it is nevertheless  possible to construct a process that satisfies \eqref{B1}, by deleting points which have a too close neighbor. As those points will be few for small volume fractions, this will not affect the discussion of \eqref{B2} qualitatively.

\mspace
We will show the following result:  for $x_1, \dots, x_n$ i.i.d. with a  law $\rho \in L^\infty$ ($\rho \geq 0$, $\int \rho = 1$),  for all $\eta > 0$:
\begin{align} \label{eq:B2.iid}
n^{-1} \# \{ i, \: \exists j \neq i, \:  |x_i - x_j| \le \eta n^{-1/3} \}  \: \xrightarrow[n \rightarrow +\infty]{} \:  1 -   \int_{\R^3} \rho(x)  e^{-\rho(x) \frac{4\pi}{3}  \eta^3} dx  
\end{align}
in probability. This implies \eqref{B2} with $\alpha = 3$ in probability. We first set
$$\eta_n := \eta n^{-1/3}, \quad   B^n_j := B(x_j, \eta_n), \quad    Y^n_{i}  :=  \prod_{j\neq i} 1_{(B^n_j)^c}(x_i).$$
Note that the random variables  $Y^n_{i}$ are identically distributed, but not independent.  Note also that $ \displaystyle n^{-1} \# \{ i, \: \exists j \neq i, \:  |x_i - x_j| \le \eta_n \} = \frac{1}{n} \sum_{i=1}^n (1-Y_i^n)$. Hence, we need to show that $\frac{1}{n} \sum_{i=1}^n Y_i^n$ converges to   $\displaystyle I_{\rho, \eta} := \int_{\R^3} \rho(x)  e^{-\rho(x) \frac{4\pi}{3}  \eta^3} dx $ in probability.

\mspace
\noindent {\em Step 1}. We show that $ \E Y_1^n  \: \xrightarrow[n \rightarrow +\infty]{} \:  I_{\rho,\eta}$. Indeed, by independence, 
\begin{align*}
  \E Y_1^n &  = \int_{\R^3} \Big(  \int_{\R^3} 1_{B(y,\eta_n)^c}(x) \rho(y)dy \Big)^{n-1} \rho(x) dx   = \int_{\R^3}   \Big( 1-  \int_{B(x, \eta_n)} \rho(y) dy \Big)^{n-1} \rho(x) dx 
  \end{align*}
 At each Lebesgue point $x$ of $\rho$, one has $\frac{1}{|B(x, \eta_n)|} \int_{B(x, \eta_n)} \rho(y) dy \rightarrow \rho(x)$, so that 
  $$   \Big( 1-  \int_{B(x, \eta_n)} \rho(y) dy \Big)^{n-1}  \rightarrow  e^{-\rho(x) \frac{4\pi}{3}  \eta^3} \quad \text{for a.e. $x$} $$
and the result follows by the dominated convergence theorem. 

\mspace
\noindent {\em Step 2}.  We show that 
$$\text{var} \Big( \frac{1}{n}\sum_{i=1}^n Y^n_i \Big) \rightarrow 0 \quad \text{as} \: n \rightarrow +\infty$$
By Markov inequality and Step 1, this implies \eqref{eq:B2.iid}.

\mspace
We have 
\begin{align*}
\text{var}\Big( \frac{1}{n} \sum_{i=1}^n Y^n_i \Big)  & = \frac{1}{n}\text{var}Y^n_1 + 
\frac{n(n-1)}{n^2} \text{Cov}Y_1^n Y_2^n = \text{Cov}Y_1^n Y_2^n + O(\frac{1}{n})
\end{align*}
using that $0 \le Y_1^n \le 1$. It remains to show that the covariance goes to zero. Using the independence of the $x_i$'s, we have the explicit formula
\begin{align*}
\E Y_1^n Y_2^n = \int_{\R^6} \Big( \int_\R 1_{|x-x_1|\ge \eta_n} 1_{|x-x_2|\ge \eta_n} \rho(x) dx \Big)^{n-2}  1_{|x_1-x_2|\ge \eta_n} \rho(x_1) \rho(x_2) dx_1 dx_2. \\
\end{align*}
We have
\begin{align*}
&\Big( \int_\R 1_{|x-x_1|\ge \eta_n} 1_{|x-x_2|\ge \eta_n} \rho(x) dx \Big)^{n-2}  \\
= & \Big( 1 -  \int_{B(x_1,\eta_n)} \!\rho -  \int_{B(x_2,\eta_n)}\!\rho + \Big(\int_{B(x_1,\eta_n) \cap B(x_2,\eta_n)} \!\! \rho\Big) \,  1_{|x_1-x_2|\le 2\eta_n}\Big)^{n-2}\\
= & e^{-\frac{4\pi}{3}\dashint_{B(x_1,\eta_n)} \rho} e^{-\frac{4\pi}{3}\dashint_{B(x_1,\eta_n)} \rho} e^{R_n(x_1,x_2)}, \quad       |R_n(x_1,x_2)|\le C  1_{|x_1-x_2|\le 2\eta_n} + C n^{-1}.
\end{align*}
This quantity converges almost surely to $e^{-\frac{4\pi}{3}\rho(x_1)} e^{-\frac{4\pi}{3}\rho(x_2)}$ and it follows by the dominated convergence theorem that 
\begin{align*}
 \E Y_1^n Y_2^n \rightarrow \Big(\int_\R e^{-\frac{4\pi}{3}\rho(x_1)} \rho(x_1) dx_1 \Big) \Big(\int_\R e^{-\frac{4\pi}{3}\rho(x_2)}\rho(x_2) dx_2 \Big) = \lim_{n \rightarrow +\infty} (\E Y_1^n)^2 
 \end{align*} 
wich yields the result.

\section*{Acknowledgements}
The first author acknowledges the support of the Institut Universitaire de France, and of the SingFlows project, grant ANR-18-CE40-0027 of the French National Research Agency (ANR)

The second author has been supported by the Deutsche Forschungsgemeinschaft (DFG, German Research Foundation) 
through the collaborative research center ``The Mathematics of Emerging Effects'' (CRC 1060, Projekt-ID 211504053) 
and the Hausdorff Center for Mathematics (GZ 2047/1, Projekt-ID 390685813).


\begin{thebibliography}{10} 


\bibitem{AlBr}
Y.~Almog and H.~Brenner.
\newblock Global homogenization of a dilute suspension of spheres. 
\newblock  Preprint arxiv:2003.01480, 2020.


\bibitem{BaGr1}
G.~Batchelor and J.~Green.
\newblock The determination of the bulk stress in a suspension of spherical
  particles at order $c^2$.
\newblock {\em J. Fluid Mech.}, 56:401--427, 1972.


\bibitem{MR1950431}
D.~J. Daley and D.~Vere-Jones.
\newblock {\em An introduction to the theory of point processes. {V}ol. {I}}.
\newblock Probability and its Applications (New York). Springer-Verlag, New
  York, second edition, 2003.
\newblock Elementary theory and methods.




\bibitem{DuerinckxGloria19}
M. Duerinckx and A. Gloria.
\newblock  Corrector equations in fluid mechanics: Effective viscosity of colloidal suspensions.
\newblock  Preprint arXiv:1909.09625, 2019.

\bibitem{DuerinckxGloria20}
M. Duerinckx and A. Gloria.
\newblock  On Einstein's effective viscosity formula.
\newblock  Preprint arXiv:2008.03837, 2020.

\bibitem{Duerinckx20}
M. Duerinckx.
\newblock  Effective viscosity of random suspensions without uniform separation.
\newblock  Preprint arXiv:2008.13188, 2020.


\bibitem{Ein}
A.~Einstein.
\newblock Eine neue Bestimmung der Molek\"uldimensionen.
\newblock {\em Ann. Physik.}, 19:289--306, 1906.

\bibitem{Galdi}
G.~Galdi.
\newblock {\em An introduction to the mathematical theory of the Navier-Stokes equations: Steady-state problems.}
\newblock Springer Science \& Business Media, 2011.

\bibitem{DGV}
D.~G\'{e}rard-Varet. 
\newblock A simple justification of effective models for conducting or fluid media with dilute spherical inclusions.
\newblock Preprint arXiv:1909.11931, 2019. 

\bibitem{Gerard-Varet20}
D.~G\'{e}rard-Varet. 
\newblock Derivation of Batchelor-Green formula for random suspensions.
\newblock Preprint arXiv:2008.06324, 2020. 


\bibitem{DGV_MH}
D.~G\'{e}rard-Varet and M.~Hillairet.
\newblock Analysis of the viscosity of dilute suspensions beyond {E}instein's
  formula.
\newblock  {\em Arch. Rat. Mech. Anal.}, 238:1349–-1411, 2020.

\bibitem{GerMec20}
D.~G\'{e}rard-Varet and A.~Mecherbet.
\newblock On the correction to Einstein's formula for the effective viscosity.
\newblock arXiv:2004.05601, 2020.



\bibitem{Guaz}
E. Guazzelli and O. Pouliquen. 
\newblock Rheology of dense granular suspensions.
\newblock {\em Journal of Fluid Mechanics}, Cambridge University Press (CUP), 2018, 852


\bibitem{MR2982744}
B.~M. Haines and A.~L. Mazzucato.
\newblock A proof of {E}instein's effective viscosity for a dilute suspension
  of spheres.
\newblock {\em SIAM J. Math. Anal.}, 44(3):2120--2145, 2012.


\bibitem{HiWu}
M.~Hillairet and D.~Wu.
\newblock Effective viscosity of a polydispersed suspension.
\newblock \{em J. Math. Pures Appl. 138:413--447, 2020.


\bibitem{Hof2}
R.~M. H\"{o}fer.
\newblock Convergence of the method of reflections for particle suspensions in Stokes flows. 
\newblock Preprint arXiv:1912.04388, 2019.  

\bibitem{HoeferSchubert20}
R.~M. H\"{o}fer and R.~Schubert.
The Influence of Einstein's Effective Viscosity on Sedimentation at Very Small Particle Volume Fraction.
\newblock Preprint arXiv:2008.04813, 2020.  



\bibitem{MR1329546}
V.~V. Jikov, S.~M. Kozlov, and O.~A. Ole\u{\i}nik.
\newblock {\em Homogenization of differential operators and integral
  functionals}.
\newblock Springer-Verlag, Berlin, 1994.
\newblock Translated from the Russian by G. A. Yosifian.


\bibitem{MR813657}
T.~L\'{e}vy and E.~S\'{a}nchez-Palencia.
\newblock Einstein-like approximation for homogenization with small
  concentration. {II}. {N}avier-{S}tokes equation.
\newblock {\em Nonlinear Anal.}, 9(11):1255--1268, 1985.


\bibitem{NiSc}
B.~Niethammer and R.~Schubert.
\newblock A local version of {E}instein's formula for the effective viscosity
  of suspensions.
\newblock {\em SIAM J. Math. Anal.}, 52(3):2561-–2591, 2020.



\bibitem{MR813656}
E.~S\'{a}nchez-Palencia.
\newblock Einstein-like approximation for homogenization with small
  concentration. {I}. {E}lliptic problems.
\newblock {\em Nonlinear Anal.}, 9(11):1243--1254, 1985.




\end{thebibliography}
\end{document}